\documentclass[12pt]{article} 

\usepackage{amsmath, amssymb, delarray, graphicx, natbib, pstricks,
pst-plot, subfigure}

\newcommand{\rcp}{\overline}

\begin{document}

\title{The Plimpton 322 Tablet and the Babylonian Method of Generating
  Pythagorean Triples} 
\author{Abdulrahman A. Abdulaziz \\ 
{\it University of Balamand} \\
}
\date{}
\maketitle

\abstract{Ever since it was published by Neugebauer and Sachs in 1945, the Old
  Babylonian tablet known as Plimpton 322 has been the subject of numerous
  studies leading to different and often conflicting interpretations of
  it. Overall, the tablet is more or less viewed as a list of fifteen
  Pythagorean triplets, but scholars are divided on how and why the list was
  devised. In this paper, we present a survey of previous attempts to
  interpret Plimpton 322, and then offer some new insights that could help in
  sharpening the endless debate about this ancient tablet.}

\section{Introduction}
Plimpton 322 is the catalog name of an Old Babylonian (OB) clay tablet held
at Columbia University. The tablet is named after New York publisher George
A.  Plimpton who purchased it from archaeology dealer Edgar J. Banks in the
early nineteen twenties.  In the mid thirties, the tablet, along with the
rest of Mr. Plimpton collection, was donated to Columbia University.
According to Banks, the tablet was found at Tell Senkereh, an archaeological
site in southern Iraq corresponding to the ancient Mesopotamian city of
Larsa [\citealt{Robson3}].

\begin{figure}[ht]
\centering
\includegraphics[width=12.7cm, height=8.8cm]{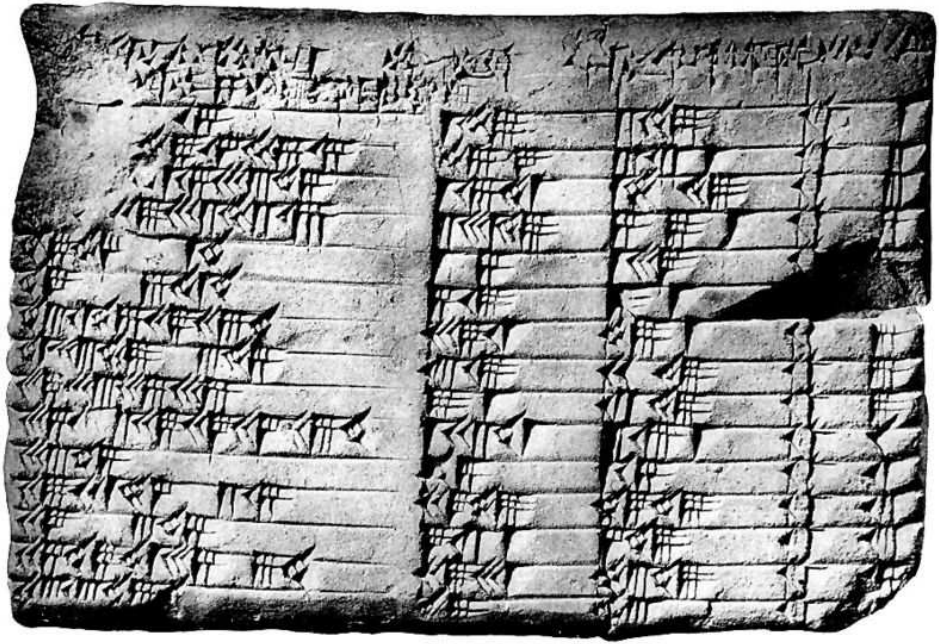}
\caption{The Plimpton 322 tablet (roughly to scale).} \label{Tablet} 
\end{figure}

The preserved portion of the tablet (shown in Figure \ref{Tablet}) is
approximately 13 cm wide, 9 cm high and 2 cm deep.
As can be seen from the picture, a chunk of the tablet is missing from the
middle of the right-hand side. Also, the tablet had (before it was baked for
preservation) remnants of modern glue on its damaged left-hand side
suggesting that it might be a part of a larger tablet, the discovery of the
remainder of which, if it ever existed, might settle many of the questions we
try to answer in this paper.  The exact date of the tablet is not known, but
it is generally agreed that it belongs to the second half of the Old
Babylonian period, roughly between 1800 and 1600 BCE. More recently, based on
the style of cuneiform script used in the tablet and comparing it with other
dated tablets from Larsa, Eleanor Robson has narrowed the date of Plimpton
322 to a period ranging from 1822 to 1784 BCE [\citealt{Robson2}].

The preserved part of Plimpton 322 is a table consisting of sixteen rows and
four columns. The first row is just a heading and the fourth (rightmost)
column of each row below the heading is simply the number of that row.  The
remaining entries are pure numbers written in \emph{sexagesimal} (base 60)
notation. However, it should be noted that due to the broken left edge of the
tablet, it is not fully clear whether or not 1 should be the leading digit of
each number in the first column. In Table \ref{Observe} we list the numbers
on the obverse of the tablet, with numbers in brackets being extrapolated. At
all times, it should be kept in mind that the original tablet would still be
of acceptable size if one or two columns were added to its left edge.
\begin{table}[ht]
\centering
\setlength{\extrarowheight}{-2pt}
\addtolength{\tabcolsep}{10pt}
\begin{tabular}{lrrr}
\hline \\[-11pt]
I& II& III& IV \\
\hline \\[-11pt]
[1 59 00] 15&                    1 59&     2 49&     1 \\[0pt]  
[1 56 56] 58 14 56 15&          56 07&   3 12 1&     2 \\[0pt] 
[1 55 07] 41 15 33 45&        1 16 41&  1 50 49&     3 \\[0pt] 
[1] 5[3] 10 29 32 52 16&      3 31 49&  5 09 01&     4 \\[0pt] 
[1] 48 54 01 40&                 1 05&     1 37&     5 \\[0pt] 
[1] 47 06 41 40&                 5 19&     8 01&   [6] \\[0pt] 
[1] 43 11 56 28 26 40&          38 11&    59 01&     7 \\[0pt] 
[1] 41 33 59 03 45&             13 19&    20 49&     8 \\[0pt] 
[1] 38 33 36 36&                 9 01&    12 49&     9 \\[0pt] 
[1] 35 10 02 28 27 24 26 40&  1 22 41&  2 16 01&    10 \\[0pt]
[1] 33 45&                         45&     1 15&    11 \\[0pt]     
[1] 29 21 54 02 15&             27 59&    48 49&    12 \\[0pt]
[1] 27 [00] 03 45&            7 12 01&     4 49&    13 \\[0pt]
[1] 25 48 51 35 06 40&          29 31&    53 49&    14 \\[0pt]
[1] 23 13 46 40&                   56&       53&  [15] \\[2pt]
\hline
\end{tabular}
\caption{The numbers on the obverse of Plimpton 322.}
\label{Observe}  
\end{table}

The numbers on the Plimpton tablet are written in cuneiform script using the
sexagesimal number system. Strictly speaking, the Babylonian number system is
not a pure sexagesimal system in the modern sense of the word. First, the
digits from 1 to 59 are expressed using only two symbols: A narrow wedge
representing 1 and a wide wedge representing 10. The numbers from 1 to 9 are
expressed by grouping the corresponding number of narrow wedges, and the
multiples of ten up to fifty are expressed by grouping the corresponding
number of wide wedges.
Every other digit is expressed as a group of wide wedges followed by a group
of narrow wedges. Second, despite the occasional indication of zero by an
empty space, the tablet lacks a consistent symbol for zero, as is typical of
mathematical texts of the same period.  Third, the Babylonian number system
does not explicitly specify the power of sixty multiplying the leading digit
of a given number. However, we will see that this uncertainty could be an
important benefit, especially since the base power can often be deduced from
context. To avoid these ambiguities as we transliterate cuneiform numbers
into modern symbols, a semicolon will be used to distinguish the whole part
from the fractional part of the number, while an empty space will be used as
a separator between the digits of the number (other authors use commas or
colons to separate the digits). We will also insert a zero wherever
necessary.

The reason behind the Babylonian use of this strange system of counting is
still debated, but one sure thing about the number sixty is that it is the
smallest number with twelve divisors: 1, 2, 3, 4, 5, 6, 10, 12, 15, 20, 30,
and 60.  With such a large number of divisors, many commonly used fractions
have simple sexagesimal representations. In fact, it is true that the
reciprocal of any number that divides a power of sixty will have a finite
sexagesimal expansion. These are the so called \emph{regular numbers}. In
modern notation, a regular number must be of the form $2^\alpha 3^\beta
5^\gamma$, where $\alpha$, $\beta$ and $\gamma$ are integers, not
necessarily positive. The advantage of 60 over $30 = 2 \times 3 \times 5$ is
that 60, unlike 30, is divisible by the highly composite number 12.

Let $m$ and $n$ be two natural numbers and denote the reciprocal of $n$ by
$\rcp n$. It follows that $\rcp n$ has a finite representation in base $60$
if and only if $n$ is regular. This is helpful because the Babylonians found
$m/n$ by computing $m$ times $\rcp n$.  To facilitate their multiplications,
they extensively used tables of reciprocals. In particular, the
\emph{standard} table of reciprocals lists the regular numbers up to 81 along
with the sexagesimal expansion of their reciprocals, as shown in Table
\ref{Frac}.
\begin{table}[ht]
\centering
\setlength{\extrarowheight}{-3pt}
\addtolength{\tabcolsep}{12pt}
\begin{tabular}{rlllllll}
\hline
$n$& $\rcp n$&& $n$& $\rcp n$&&     $n$&   $\rcp n$ \\
\hline
 2&        30&&  16&     3 45&&      45&     \ 1 20 \\ 
 3&        20&&  18&     3 20&&      48&     \ 1 15 \\
 4&        15&&  20&        3&&      50&     \ 1 12 \\      
 5&        12&&  24&     2 30&&      54&  \ 1 06 40 \\     
 6&        10&&  25&     2 24&&     \ 1&        \ 1 \\
 8&    \ 7 30&&  27&  2 13 20&&  \ 1 04&      56 15 \\
 9&    \ 6 40&&  30&        2&&  \ 1 12&         50 \\
10&       \ 6&&  32&  1 52 30&&  \ 1 15&         48 \\
12&       \ 5&&  36&     1 40&&  \ 1 20&         45 \\
15&       \ 4&&  40&     1 30&&  \ 1 21&   44 26 40 \\
\hline
\end{tabular}
\caption{The standard Babylonian table of reciprocals.} \label{Frac}
\end{table}
Fractions like $\rcp7$ are omitted from the table because they do not have
finite expansions, but sometimes approximations were used for such small
non-regular numbers. For example,
\[
\rcp7 = 13 \times \rcp{91}  \approx 13 \times
\rcp{90} = 13 \times (0;00 \ 40) = 0;08 \ 40.
\]
The Babylonians went even further than this. On another OB tablet published
by A. Sachs, we find a lower and an upper bound for $\rcp7$
[\citealt{Sachs}]. It states what we now write as
\[
0;08 \ 34 \ 16 \ 59 < \rcp7 < 0;08 \ 34 \ 18. 
\]
Amazingly, the correct value of $\rcp7$ is $0;08 \ 34 \ 17$, where the part to
the right of the semicolon is repeated indefinitely.

\section{Interpretation of the  Tablet} 
Originally, Plimpton 322 was classified as a record of commercial
transactions. However, after Neugebauer and Sachs gave a seemingly
irrefutable interpretation of it as something related to Pythagorean
triplets, the tablet gained so much attention that it has probably become the
most celebrated Babylonian mathematical artifact [\citealt{Neugebauer1}]. For
some, like Zeeman, the tablet is hailed as an ancient document on number
theory, while for others, like Robson, it is just a school record of a
student working on selected exercises related to squares and reciprocals
[\citealt{Zeeman}; \citealt{Robson2}]. But no matter which view one takes,
there is no doubt that the tablet is one of the greatest achievements of OB
mathematics, especially since we know that it was written at least one
thousand year before Pythagoras was even born.  That the Old Babylonians knew
of Pythagoras theorem (better called rule of right triangle) is evident in
the many examples of its use in various problems of the same period
[\citealt{Hoyrup2}].  Having said this, it should be clear that the tablet is
no way a proof of Pythagoras theorem. In fact, the idea of a formal proof is
nowhere to be found in extant Babylonian mathematics [\citealt{Friberg1}].

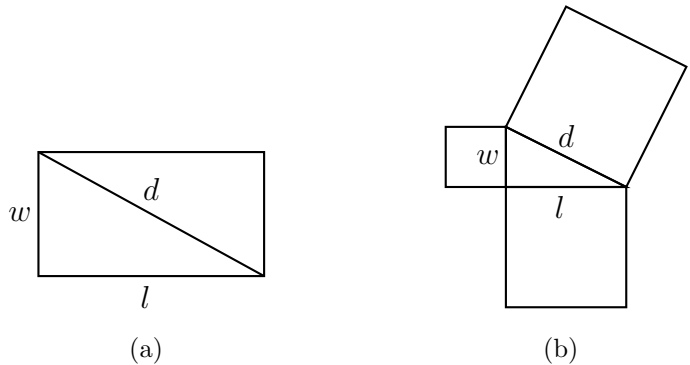
\begin{figure}[t]
\centering
\psset{xunit=0.75cm}
\psset{yunit=0.55cm}
\subfigure[]{
\begin{pspicture}(0,-0.75)(4,3.25)
\pspolygon(0,0)(4,0)(4,3)(0,3)
\psline(0,3)(4,0)
\rput(1.9,-0.5){$l$}
\rput(-0.33,1.5){$w$}
\rput(2,2.1){$d$}
\end{pspicture}
}
\hskip 2cm
\psset{xunit=0.8cm}
\psset{yunit=0.8cm}
\subfigure[]{
\begin{pspicture}(0,-2)(4,3)
\pspolygon(1,0)(3,0)(1,1)
\pspolygon(0,0)(1,0)(1,1)(0,1)
\pspolygon(1,0)(3,0)(3,-2)(1,-2)
\pspolygon(1,1)(2,3)(4,2)(3,0)
\rput(0.7,0.5){$w$}
\rput(1.9,-0.33){$l$}
\rput(2,0.85){$d$}
\end{pspicture}
}
\caption{(a) The diagonal as the hypotenuse and (b) the sides as squares.}
\label{RectPyth} 
\end{figure}

According to Neugebauer and Sachs, the heading of the fourth column is `its
name', which simply indicates the line number, from 1 to 15
[\citealt{Neugebauer1}].  The headings of columns two and three read
something like `square of the width (or short side)' and `square of the
diagonal', respectively. An equally consistent interpretation can be obtained
if the word `square' is replaced by `square root'. These headings make sense
only when coupled with the fact that the Babylonian thought of the sides of a
right triangle as the length and width of a rectangle whose diagonal is the
hypotenuse of the given triangle. Also, the Babylonians used the word square
to mean the side of a square as well as the square itself
[\citealt{Robson2}]. Let the width, length and diagonal of the rectangle be
denoted by $w$, $l$ and $d$. Then the relation between the right triangle and
the rectangle is shown in Figure \ref{RectPyth}(a), while in Figure
\ref{RectPyth}(b) three squares are drawn, one for each side of the triangle.
Indeed, Figure \ref{RectPyth}(b) should look familiar to anyone acquainted
with Euclid's proof of Pythagoras theorem. Such Babylonian influence on Greek
mathematics is in accordance with tales that Pythagoras spent more than
twenty years of his life acquiring knowledge from the wise men of Egypt and
Mesopotamia [\citealt[p.~85]{Bell}].

So, using Neugebauer's interpretation, the second column represents the short
side of a right triangle or the width $w$ of the corresponding rectangle and
the third column represents the hypotenuse of the right triangle or the
diagonal $d$ of the rectangle. The longer side of the triangle or length $l$
of the rectangle does not appear in the table (maybe it was written on the
missing part of the tablet). In such an interpretation, the first column is
simply $d^2/l^2$ or $(d/l)^2$. To see how this can be applied to the table,
let us look at an example. In line 5, the square of the number in Column III
minus the square of the number in Column II is a perfect square. That is,
\[
(1\ 37)^2 - (1\ 05)^2 = 1\ 26\ 24 = (1\ 12)^2.
\]  
Moreover, the number in the first column is nothing but the square of the
ratio (1 37):(1 12). In decimal notation, we have
\[
97^2 - 65^2 = 5184 = 72^2, 
\]
with the number in Column I being $97^2/72^2$. 

In general, if we think of the two middle entries in a given row as the width
$w$ and diagonal $d$ of a rectangle (or the short side and hypotenuse of a
right triangle), then the entry in the first column of that row is $d^2/l^2$,
where $l$ is the length of the rectangle. Except for a few errors, which we
will say more about in Section \ref{ErrorSec}, the entries in the first three
columns of each line of Table \ref{Observe} are exactly $d^2/l^2$, $d$ and
$w$.  It is important to note that an equally consistent interpretation can
be obtained if we remove the leading 1 in Column I and think of the numbers
in that column as $w^2/l^2$.\footnote{After personally inspecting Plimpton
  322, Bruins concluded that the apparent unit at the beginning of each line
  is due to the horizontal line between rows [\citealt{Bruins1}].  However,
  Friberg has given a more convincing argument in support of the leading one
  stance [\citealt{Friberg1}].} The equivalence of the two interpretations
follows from the fact that if $l^2 + w^2 = d^2$, then
\begin{equation}
1+w^2/l^2 = d^2/l^2. \label{ColI}
\end{equation}

This brings us to the heading of the first column. The heading consists of
two lines, both of which are damaged at the beginning. Quoting Neugebauer and
Sachs:
\begin{quote}
The translation causes serious
  difficulties. The most plausible rendering seems to be: `The
  \emph{takiltum} of the diagonal which has been subtracted such that the
  width. \dots'  [\citealt[p. 40]{Neugebauer1}]
\end{quote}
Building on the work of H\o yrup, Robson took the work of Neugebauer and
Sachs a step further [\citealt{Hoyrup1}; \citealt{Robson2}]. By carefully
examining the damaged heading of the first column, she was able to render the
sensible translation:
\begin{quote}
The holding-square of the diagonal from which 1 is torn out, so that the
short side comes up.
\end{quote}
In addition to being linguistically and contextually sound, the above
interpretation, thought of as equation (\ref{ColI}) in words, makes perfect
mathematical sense.\footnote{Price mentioned a similar interpretation proposed
  by Goetze [\citealt{Price}].} This is in line with the `cut-and-paste'
geometry introduced by H\o yrup and adopted by Robson.  Also, the
interpretation clearly speaks in favor of the restoration of the leading 1 at
the beginning of each line of the preserved tablet.

The above interpretation of the tablet is the most widely accepted one by
scholars because it relates the numbers on the tablet in a meaningful way
which is totally drawn from extant OB mathematics. Other, wilder,
interpretations of the tablet have also been proposed but they do not carry
much weight. One such interpretation suggests that the tablet represents some
sort of a trigonometric table [\citealt{Joyce}; \citealt[pp.
30-34]{Maor}]. This is based on the fact that in each line of the Plimpton
table, the entry in the first column is the square of the cosecant of the
angle between the long side and the hypotenuse of a right triangle, where the
angle decreases from about $45^\circ$ to $30^\circ$ by roughly one degree per
line as we move down the table. However, this hypothesis is dismissed by most
historians of Babylonian mathematics on many grounds, the least of which is
the lack of any traces of trigonometric functions in extant Babylonian
mathematics. Moreover, the translation of the heading of the first column
offered by Robson totally refutes such an interpretation [\citealt{Robson2}].

\section{Previous Methods for Reconstructing the Table}
Having determined what the tablet means, it remains to answer the more
difficult question of how it was constructed. In this respect, there are two
major theories on how the numbers on Plimpton~322 were generated. The first
method was proposed by Neugebauer and Sachs in their highly acclaimed book
\emph{Mathematical Cuneiform Texts} in which the tablet was originally
published [\citealt[pp.~38-~41]{Neugebauer1}]. The method had many proponents
including [\citealt{Gillings1}; \citealt{Price}; \citealt{Buck}] and
others. The second method was first introduced by E.  M. Bruins in 1949, but
it did not become main stream until it reappeared in the works of Schmidt and
Friberg, and more recently in the work of Robson [\citealt{Schmidt};
\citealt{Friberg1}; \citealt{Robson2}]. The decision of which method was
employed in the construction of the tablet is made more difficult by the fact
that the two methods are mathematically equivalent to each other.  As a
general rule, one should pick the method which is more consistent with extant
mathematics of the OB period. But before this could be done, we shall give
some background information and a summary of each method.

It was known to the ancient Greeks that if $m > 1$ is an odd integer, then
the numbers
\begin{equation*}
m, \qquad \frac{1}{2}(m^2-1) \qquad \text{and} \qquad \frac{1}{2}(m^2+1)
\end{equation*}
satisfy Pythagoras theorem [\citealt[p. 121]{Shanks}]. That is,
\begin{equation*}
m^2 + \left[\frac{(m^2-1)}{2}\right]^2 =
\left[\frac{(m^2+1)}{2}\right]^2. 
\end{equation*}
The restriction that $m$ is an odd integer is needed to guarantee that the
triplet consists of integers, but the theorem holds for every real number
$\alpha$. In particular, if we replace $m$ by $\alpha > 0$ and divide by
$\alpha^2$, we obtain the normalized equation
\begin{equation}
1+\left( \frac{\alpha - 1/\alpha}{2} \right)^2 = \left( \frac{\alpha +
1/\alpha}{2} \right)^2  
\cdot \label{BaThm1}  
\end{equation}
Observe that as $\alpha$ varies between $1$ and $ 1 +\sqrt2$, $(\alpha -
1/\alpha)/2$ varies between $0$ and $1$; and when $\alpha$ increases beyond
$1 + \sqrt2$, $(\alpha - 1/\alpha)/2$ increases beyond 1. In the former case,
the longer side of the triangle must be 1; while in the latter the reverse is
true. The case $\alpha < 1$ will be ignored since it leads to negative values
of $\alpha -1/\alpha$.

Although (\ref{BaThm1}) holds for every positive real $\alpha$, we are mainly
interested in the case when $\alpha$ is a regular number. But if $\alpha$ is
regular, then it must be of the form $p/q$, where $p$ and $q$ are regular
integers.  This yields
\begin{equation}
1 + \left( \frac{p^2-q^2}{2pq} \right)^2 = \left( \frac{p^2+q^2}{2pq}
\right)^2 \cdot \label{BaThm2}
\end{equation}
Now multiplying by $(2pq)^2$, we see that (\ref{BaThm2}) is the same as
saying that the integers
\begin{equation}
2pq, \qquad p^2-q^2 \qquad \text{and} \qquad p^2+q^2
\label{Triplet2}
\end{equation}
form an integral Pythagorean triplet, a fact also known to Euclid. More
generally, if $p$ and $q$ are integers of opposite parity with no common
divisor and $p > q$, then all \emph{primitive} Pythagorean triplets are
generated by (\ref{Triplet2}) [\citealt[p. 141]{Shanks}]. By a primitive
triplet we mean a triplet consisting of numbers that are prime to each other.

There is direct evidence that the Babylonians were aware of and even used an
identity the like of (\ref{BaThm1}) or (\ref{BaThm2})
[\citealt{Bruins1}]. Also, according to Neugebauer and Sachs, the numbers in
the first three columns of the tablet are merely $d^2/l^2$, $w$ and $d$,
where $l = 2pq$, $w = p^2 - q^2$ and $d = p^2 + q^2$ are as in
(\ref{Triplet2}).
\begin{table}[ht]
\centering
\setlength{\extrarowheight}{-3pt}
\addtolength{\tabcolsep}{7pt}
\begin{tabular}{rrrlrrr}
\hline \\[-12pt]  
$p$& $q$&   $l$&              $d^2/l^2$&            $w$&            $d$& $n$\\
\hline
  12&   5&    2 00&                1 59 00 15&\hskip6pt 1 59&        \ 2 49&   1\\
1 04&  27&   57 36&    1 56 56 58 14 50 06 15&         56 07&       1 20 25&   2\\
1 15&  32& 1 20 00&       1 55 07 41 15 33 45&       1 16 41&       1 50 49&   3\\
2 05&  54& 3 45 00&       1 53 10 29 32 52 16&       3 31 49&       5 09 01&   4\\
   9&   4&    1 12&             1 48 54 01 40&          1 05&          1 37&   5\\
  20&   9&    6 00&             1 47 06 41 40&          5 19&\hskip6pt 8 01&   6\\
  54&  25&   45 00&       1 43 11 56 28 26 40&         38 11&         59 01&   7\\
  32&  15&   16 00&       1 41 33 45 14 03 45&         13 19&         20 49&   8\\
  25&  12&   10 00&             1 38 33 36 36&          8 01&         12 49&   9\\
1 21&  40& 1 48 00& 1 35 10 02 28 27 24 26 40&       1 22 41&       2 16 01&  10\\
   2&   1&       4&                   1 33 45&             3&             5&  11\\
  48&  25&   40 00&          1 29 21 54 02 15&         27 59&         48 49&  12\\
  15&   8&    4 00&             1 27 00 03 45&          2 41&          4 49&  13\\
  50&  27&   45 00&       1 25 48 51 35 06 40&         29 31&         53 49&  14\\
   9&   5&    1 30&             1 23 13 46 40&            56&          1 46&  15\\
\hline
\end{tabular}
\caption{Plimpton 322 preceded by the generators $p$ and $q$ and the long side
 $l$.}
 \label{Complete}  
\end{table}
In Table \ref{Complete}, we list the numbers on the original tablet
(corrected when necessary) along with the missing side $l$ as well as the
generating parameters $p$ and $q$.  What makes this interpretation
particularly attractive is that except for $p=2$ 05 in row four of the
table, all values of $p$ and $q$ are in the standard table of reciprocals
(Table \ref{Frac}).  But it happens that the regular number 2 05 is the first
(restored) value in another table of reciprocals found in the OB tablet CBS
29.13.21 [\citealt[p.  14]{Neugebauer1}]. Moreover, the value of $p/q$,
like that of $d^2/l^2$, decreases as we move down the table. In fact, if we
allow $p$ and $q$ to be regular numbers less than or equal to 2 05 and sort
the resulting table by $p/q$ (or by $d^2/l^2$), then only one extra triplet
that goes between lines 11 and 12 is produced, see Table~\ref{Extra}.
\begin{table}[ht]
\centering
\addtolength{\tabcolsep}{5pt}
\begin{tabular}{lllllll}
\hline
 $p$&  $q$&     $l$&              $d^2/l^2$&     $w$&     $d$& $n$ \\
\hline
2 05& 1 04& 4 26 40& 1 31 09 09 25 42 02 15& 3 12 09& 5 28 41& 11a \\
\hline
\end{tabular}
\caption{The missing line 11a should be inserted between lines 11 and 12.}
\label{Extra}  
\end{table}
This led D. E. Joyce to argue that the extra triplet may have been
inadvertently left out or that it was dismissed because the magnitudes of the
sides in the resulting triplet are too large [\citealt{Joyce}]. On the other
hand, if we consider all regular integers $p$ and $q$ such that $p \leq 2$ 05
and $q \leq 1$ 00, then the first fifteen triplets, sorted by descending
value of $p/q$, are exactly those found in the Plimpton tablet
[\citealt{Price}].

The second widely accepted method for devising the tablet was proposed by E.
M. Bruins in 1949, and is often called the `reciprocal' method as opposed
to the `generating pair' method of Neugebauer and Sachs. The method is based
on the fact that if $r = p \rcp q$ is regular, then the numbers
$x=\rcp2(r-\rcp r)$ and $y=\rcp2(r+\rcp r)$ satisfy, on top of being regular,
the equation $1+x^2 = y^2$. It turned out that for each line in the tablet
one can find a regular number $r$ (called the \emph{generating ratio}) such
that when $x$ and $y$ are divided by their regular common factors, we end up
with the numbers in the second and third columns of the Plimpton tablet. To
see how this could be done, let us take a closer look at line five of Table
\ref{Complete}.  Since $p = 9$ and $q=4$, we have $r=2;15$ and $\rcp r =
0;26\ 40$. This yields $x=0;54 \ 10$ and $y=1;20\ 50$. Since the rightmost
(sexagesimal) digit of $x$ and that of $y$ are divisible by 10, we should
multiply both $x$ and $y$ by 6 (or equivalently divide $x$ and $y$ by
10). The resulting numbers are $5;25$ and $8;05$, both of which should be
multiplied by 12 (or divided by 5) leading to 1 05 and 1 37. In tabular form,
we have
\begin{table*}[ht]
\centering
\begin{tabular}{rrllllr}
 $6$& $\times$& 0;54 10& \qquad 1;20 50& $\times$& $6$  \\
$12$& $\times$&    5;25& \qquad    8;05& $\times$& $12$ \\
    &         &    1 05& \qquad    1 37&                \\
\end{tabular}
\end{table*}

\noindent
Since the terminating digits of 1 05 and 1 37 have nothing in common, the
process stops here. Now if we think of $1 \ 05$ and $1 \ 37$ as the width $w$
and diagonal $d$ of a rectangle, then the length $l$ must be $1 \ 12$, which
is the product of 12 and 6. Equivalently, $w = p^2-q^2, \ d = p^2+q^2,$ and
$l = 2pq$. In Table \ref{Recp}, we list the values of $r$, $\rcp r$, $x$, $y$
and $l$.
\begin{table}[ht]
\centering
\setlength{\extrarowheight}{-3pt}
\addtolength{\tabcolsep}{2pt}
\begin{tabular}{rrllllrr}
\hline 
 $p$& $q$&        $r$&    $\rcp r$&         $x$&           $y$&    $l$& $n$\\
\hline                                                 
  12&  5&       2;24& 0;25         &       0;59 30&       1;24 30&    2 00& 1\\
1 04& 27& 2;22 13 20& 0;25 18 45   & 0;58 27 17 30& 1;23 46 02 30&   57 36& 2\\
1 15& 32& 2;20 37 30& 0;25 36      &    0;57 30 45&    1;23 06 45& 1 20 00& 3\\
2 05& 54& 2;18 53 20& 0;25 55 12   &    0;56 29 04&    1;22 24 16& 3 45 00& 4\\
   9&  4&       2;15& 0;26 40      &       0;54 10&       1;20 50&    1 12& 5\\
  20&  9&    2;13 20& 0;27         &       0;53 10&       1;20 10&    6 00& 6\\
  54& 25&    2;09 36& 0;27 46 40   &    0;50 54 40&    1;18 41 20&   45 00& 7\\
  32& 15&       2;08& 0;28 07 30   &    0;49 56 15&    1;18 03 45&   16 00& 8\\
  25& 12&       2;05& 0;28 48      &       0;48 06&       1;16 54&   10 00& 9\\
1 21& 40&    2;01 30& 0;29 37 46 40& 0;45 56 06 40& 1;15 33 53 20& 1 48 00&10\\
   2&  1&          2& 0;30         &          0;45&         1;15 &       4&11\\
  48& 25&    1;55 12& 0;31 15      &    0;41 58 30&    1;13 13 30&   40 00&12\\
  15&  8&    1;52 30& 0;32         &       0;40 15&       1;12 15&    4 00&13\\
  50& 27& 1;51 06 40& 0;32 24      &    0;39 21 20&    1;11 45 20&   45 00&14\\
   9&  5&       1;48& 0;33 20      &       0;37 20&       1;10 40&      45&15\\
\hline
\end{tabular}
\caption{Reciprocal method: $r = p \rcp q$, $x=\rcp2(r-\rcp r)$, 
  $y=\rcp2(r+\rcp r)$ and $l=2pq$.} \label{Recp}
\end{table}
Observe that even though $p$ and $q$ change erratically from one line to the
next, the ratio $p/q$ steadily decreases as we move down the table.  Since it
is a common Babylonian practice to list numbers in descending or ascending
order (the ratios $d^2/l^2$ in Column I of the tablet decrease from top to
bottom), we have the first piece of evidence showing that the $r$-method is
more in accordance with OB mathematics.

Although the reciprocal method (or $r$-method) may seem awkward from a modern
point of view, there is ample evidence that the techniques it employs have
been used in OB mathematics [\citealt{Robson2}; \citealt{Friberg1};
\citealt{Bruins2}]. One advantage of the method is that the concept of
relatively prime numbers becomes unnecessary: Start with any regular number
$r$, the process terminates with a primitive triplet.  Moreover, the method
is favored because it provides a simple way to compute the ratio $d^2/l^2$ in
the first preserved column. Since $d/l$ is just $y$, all the scribe had to do
was calculate $\rcp2(r+\rcp r)$ and then square the result. The method also
makes use of the flexibility inherent in the (ambiguous) Babylonian number
system, where multiplication and division can be interchanged at will.  This
is due to the dismissal of leading and trailing zeros as well as the lack of
a symbol that separates the fractional from the whole part of a number.  In
short, the advantage of Bruins' method over other methods can be summarized
as follows:
\begin{quote}
So what does make Bruins' reciprocal theory more convincing than the standard
$p$,$q$ generating function---or, indeed, the trigonometric table? I have
already showed that its starting points (reciprocal pairs, cut-and-paste
algebra) and arithmetical tools (adding, subtracting, halving, finding square
sides) are all central concerns of Old Babylonian mathematics: it is
sensitive to the ancient thought-processes and conventions in a way that no
other has even tried to be. For example, in this theory the values in Column
I are a necessary step towards calculating those in Column III and may also
be used for Column II. And the Column I values themselves are derived from an
ordered list of numbers [\citealt{Robson2}].
\end{quote}

\section{Possible Ways to Complete the Table}
Suppose that one wants to list all Pythagorean triplets $(w,l,d)$ such that
$w < l < d < 20000$, $l < 15000$ is a regular, $\gcd(l,d) = 1$ and $d^2/l^2
< 2$. Then the first $15$ triplets listed in descending order of $d^2/l^2$
are exactly those found in Plimpton 322. The only exception is that the
triplet (45, 60, 75) in the tablet should be replaced by the equivalent
triplet (3, 4, 5), since the greatest common divisor of $l$ and $d$ in the
former triplet is different from one.

Although the above triplets agree with those in Plimpton 322, no one would be
imprudent enough to think that the Babylonians constructed the tablet by
following such a predefined set of modern rules, not to mention the enormous
number of calculations involved. From the outset, it should be made
absolutely clear that it is not enough to provide a set rules that produce
the numbers in the tablet unless those rules can be found, at least
implicitly, in the collective body of OB mathematics. So in order to
determine how the tablet may have been devised, we must only work with the
sort of mathematics that OB scribes had at their disposal.

\begin{figure}[ht]
\centering
\mbox{
\includegraphics[width=4.5cm, height=4.5cm]{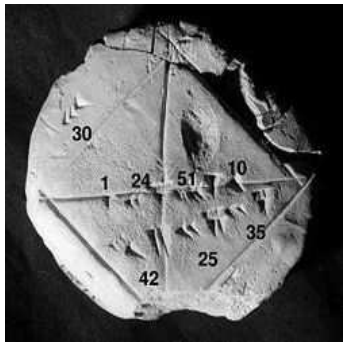}
}
\caption{The YBC 7289 tablet. Courtesy: Bill Casselman.} \label{Ybc}
\end{figure}

The Old Babylonians used 1;25 as a rough estimation of root two, and there is
strong evidence that they also used the much closer approximation of
1;24~51~10. The evidence for the closer approximation comes from two OB
tablets known as YBC 7243 and YBC 7289
[\citealt[pp. 42-43]{Neugebauer1}]. The first of the two tablets contains a
list of coefficients, where on the tenth line appears the number 1~24~51~10
followed by the words `Diagonal, square root'. The second tablet is a round
tablet circumscribing a (diagonal) square with the cuneiform symbol for 30
written above the middle of the upper-left side. In addition, the sexagesimal
numbers 1~24~51~10 and 42 25 35 are inscribed along the horizontal diagonal
and across the lower half of the vertical diagonal, see Figure~\ref{Ybc}.
Since the Old Babylonians did not explicitly write trailing zeros and since
the result of multiplying 30 by 1~24~51~10 is 42~25~35~00, it is immediately
clear that 42~25~35 should be interpreted as the product of 30 and
1~24~51~10. A more meaningful relation between the three numbers can be
deduced if we think of some or all of them as fractions and not integers.  We
can always do this because the Babylonian number system does not distinguish
between fractions and whole numbers. Taking 1~24~51~10 as 1;24~51~10 and
comparing it with $\sqrt2 = 1;24\ 51\ 10\ 07\dots$, there is little doubt
that the number at hand is a Babylonian approximation of root two.\footnote{In
  decimal notation, we have $1;24\ 51\ 10 = 1.41421\overline{296}$, while
  $\sqrt2 = 1.414213\dots$.} This makes perfect sense since the length of the
diagonal of a square is equal to radical two times the length of its
side. Furthermore, if we think of $30$ as $0;30 = \rcp 2$, then the length of
the diagonal will be equal to the reciprocal of root two. In other words, if
the length of the side is $\rcp 2$, then the length of the diagonal is a half
times root two, which is the same as the reciprocal of root two. From this we
see that the side of the square was so cleverly chosen so that the other two
numbers inscribed on the tablet are nothing but highly accurate
approximations of root two and its reciprocal. Knowing the central role
reciprocals played in Babylonian mathematics, it is hard to believe that this
could have happened by accident [\citealt{Melville2}]. As to how the
Babylonians may have found such an extremely good approximation of root two,
various ways have been proposed by different authors
[\citealt[p. 43]{Neugebauer1}; \citealt{Fowler}].

The YBC 7289 tablet is another attestation of the Old Babylonian
understanding of the theorem of Pythagoras.  The tablet covers a special case
of the theorem where the rectangle is replaced by a square, meaning that the
width $w$ is the same as the length $l$. But if $w = l$, then the equation
$w^2 + l^2 = d^2$ can be rewritten as $d^2/l^2 = 2$. The other special case
of the theorem is the extreme case with $d=l$, where the rectangle collapses
into a line. It is between these two special cases that the numbers in the
first column of the Plimpton tablet should be viewed. Let $\alpha_0$ be the
generating ratio for the first special case and, for the lack of a better
term, $\alpha_\infty$ be the generating ratio for the other case. It follows
from (\ref{BaThm1}) that $\alpha_0 = 1+\sqrt{2}$, $\alpha_\infty = 1$, and
that every positive ratio $\alpha$ such that $\alpha- 1/\alpha > 0$ must
satisfy the inequality
\begin{equation}
\alpha_0 > \alpha > \alpha_\infty, \label{Admis}
\end{equation}
provided that the length of the rectangle is taken as 1. Moreover, if
$\alpha$ is replaced by the regular number $r = p\rcp q$, then we can use
(\ref{BaThm2}) to show that for $r$ to satisfy (\ref{Admis}) we must have
\begin{equation}
p^2-q^2 < 2pq. \label{Cond}
\end{equation}
A regular ratio $r$ satisfying (\ref{Cond}) will henceforth be called
\emph{admissible}.
 
If we look at the fifteen pairs $(p,q)$ that generate the numbers inscribed
on the Plimpton tablet, we find that their ratios $p \rcp q$ are all
admissible. Moreover, the largest value of $p$ is 2 05 and that of $q$ is 54.
Assuming for now that no ratio is allowed to have larger values of $p$ and
$q$, we get a total of 38 admissible ratios, the first fifteen of which are
exactly those of Table~\ref{Recp}. The same set of ratios is produced if $q$
is allowed to be as large as sixty. In Table~\ref{Cont}, we list the $38$
ratios along with the Pythagorean triplets they generate if the $r$-method is
used. A similar table was first devised by Price using the $pq$-method
[\citealt{Price}].\footnote{Price made many calculation errors in the complete
  table.  In fact, the value of $\left(\frac{p^2+q^2}{2pq}\right)^2$ in his
  table ($d^2/l^2$ in our table) is incorrect in lines 16, 17, 23, 24, 25,
  28, 29, 30 and 34.} The lines ending with an asterisk (lines 11, 15, 18
and 36) are those for which the $pq$-method deviates from the
$r$-method. These lines correspond to ratios where $p$ and $q$ are both odd,
and consequently $p^2 - q^2$ and $p^2 + q^2$ are both even.  It follows that
$p^2-q^2$, $p^2+q^2$ and $2pq$ have a common factor of $2$ and so the
Pythagorean triplet is not primitive.  In such cases, the values produced by
the $pq$-method will be twice the values produced by the $r$-method.

\begin{table}[ht]
\centering
\setlength{\extrarowheight}{-5pt}
\addtolength{\tabcolsep}{5pt}
\begin{tabular}{lrlrrl}
\hline \\[-10pt] 
       $r$&     $l$&                 $d^2/l^2$&     $w$&      $d$& $n$  \\[2pt]
\hline \\[-12pt] 
      2;24&    2 00&                1;59 00 15&    1 59&     2 49& 1   \\
2;22 13 20&   57 36&    1;56 56 58 14 50 06 15&   56 07&  1 20 25& 2   \\
2;20 37 30& 1 20 00&       1;55 07 41 15 33 45& 1 16 41&  1 50 49& 3   \\
2;18 53 20& 3 45 00&       1;53 10 29 32 52 16& 3 31 49&  5 09 01& 4   \\
      2;15&    1 12&             1;48 54 01 40&    1 05&     1 37& 5   \\
   2;13 20&    6 00&             1;47 06 41 40&    5 19&     8 01& 6   \\
   2;09 36&   45 00&       1;43 11 56 28 26 40&   38 11&    59 01& 7   \\
      2;08&   16 00&       1;41 33 45 14 03 45&   13 19&    20 49& 8   \\
      2;05&   10 00&             1;38 33 36 36&    8 01&    12 49& 9   \\
   2;01 30& 1 48 00& 1;35 10 02 28 27 24 26 40& 1 22 41&  2 16 01& 10  \\
         2&    1 00&                   1;33 45&      45&     1 15& 11* \\
   1;55 12&   40 00&          1;29 21 54 02 15&   27 59&    48 49& 12  \\
   1;52 30&    4 00&             1;27 00 03 45&    2 41&     4 49& 13  \\
1;51 06 40&   45 00&       1;25 48 51 35 06 40&   29 31&    53 49& 14  \\
      1;48&      45&             1;23 13 46 40&      28&       53& 15* \\
\hline \\[-8pt] 
   1;46 40&    4 48&          1;22 09 12 36 15&    2 55&     5 37& 16  \\
   1;41 15&   14 24&       1;17 58 56 24 01 40&    7 53&    16 25& 17  \\
      1;40&      15&                   1;17 04&       8&       17& 18* \\
   1;37 12& 2 15 00& 1;15 04 53 43 54 04 26 40& 1 07 41&  2 31 01& 19  \\
      1;36&    1 20&             1;14 15 33 45&      39&     1 29& 20  \\
   1;33 45&   13 20&          1;12 45 54 20 15&    6 09&    14 41& 21  \\
      1;30&      12&                   1;10 25&       5&       13& 22  \\
1;28 53 20&   36 00&       1;09 45 22 16 06 40&   14 31&    38 49& 23  \\
   1;26 24&   30 00&             1;08 20 16 04&   11 11&    32 01& 24  \\
   1;25 20& 1 36 00&    1;07 45 23 26 38 26 15&   34 31&  1 42 01& 25  \\
1;24 22 30&   48 00&       1;07 14 53 46 33 45&   16 41&    50 49& 26  \\
   1;23 20&   15 00&             1;06 42 40 16&    5 01&    15 49& 27  \\
      1;21&   18 00&       1;05 34 04 37 46 40&    5 29&    18 49& 28  \\
      1;20&      24&                1;05 06 15&       7&       25& 29  \\
   1;16 48&   26 40&       1;03 43 52 35 03 45&    6 39&    27 29& 30  \\
      1;15&      40&                1;03 02 15&       9&       41& 31  \\
      1;12&    1 00&                   1;02 01&      11&     1 01& 32  \\
1;11 06 40&   28 48&       1;01 44 55 12 40 25&    4 55&    29 13& 33  \\
   1;07 30&    2 24&             1;00 50 10 25&      17&     2 25& 34  \\
   1;06 40&    3 00&             1;00 40 06 40&      19&     3 01& 35  \\
   1;04 48&   11 15&       1;00 21 21 53 46 40&      52&    11 17& 36* \\
      1;04&    8 00&          1;00 15 00 56 15&      31&     8 01& 37  \\
   1;02 30&   20 00&             1;00 06 00 09&      49&    20 01& 38  \\
\hline
\end{tabular}
\caption{The continuation of Plimpton 322. The first two columns should be
  inscribed on the broken part of the tablet; while lines 16 to 38 should be
  inscribed on the reverse of the tablet.} \label{Cont} 
\end{table}
\clearpage

To better understand the difference between the two methods, let us take a
closer look at lines 11 and 15. For line 11, Price took $p=1\ 00$ and $q=30$
so that $p\rcp q = 2$, yielding the triplet (1 00 00, 45 00, 1 15 00). This
seems somewhat contrived since taking $p=2$ and $q=1$ produces the equivalent
triplet (4, 3, 5). On the other hand, Bruins method yields $x=0;45$ and $y =
1;15$. Since these two numbers along with $l=1$ make the well known triplet
(1 00,~45,~1~15) the scribe did not bother do the simplification to obtain
the reduced triplet (4, 3, 5).  Although it is hard to decide which method
was used based only on this case, we still think that the $r$-method better
explains why the non-reduced triplet is the one that eventually appeared on
the tablet. We take this view in light of the fact that the non-reduced
triplet (1, 45, 1 15) can be found in another OB text from Tell Dhiha'i
[\citealt{Baqir}]. The text poses and solves the following problem: Find the
sides of the rectangle whose diagonal is 45 and whose area is 1 15. The
relevance of the problem to our case lies not only in the fact that the
calculated sides and the diagonal form the triplet (1,~45, 1~15), but more
importantly in the clear resemblance between the solution algorithm and the
$r$-method [\citealt{Friberg1}]. As for line 15, the $pq$-method yields the
triplet (1~30, 56, 1 46), while the $r$-method yields the triplet (45, 28,
53). Since the part of the triplet inscribed on the tablet is $w=56$ and
$d=53$, it is not clear which entry should be considered as the wrong one.
If the $pq$-method is used then the erroneous value would be that of $d$;
while if the $r$-method is used then the incorrect value would be that of $w$
(more on this in Section \ref{ErrorSec}). Had the scribe continued to line 18
or 36, we would have been able to tell which method was used with a higher
degree of certainty.  Unfortunately, lines 18 and 36 are not inscribed,
forcing us to ponder over the criteria by which the number of lines in the
tablet was determined.

If there are 38 admissible ratios with $p \leq 2 \ 05$ and $q \leq 1 \ 00$,
then one wonders why only 15 of those are found on Plimpton 322.  For
proponents of the trigonometric table it is because these ratios roughly
correspond to angles between $45^\circ$ and $30^\circ$; while for proponents
of the incomplete table the missing ratios should have also been inscribed,
perhaps on the reverse of the tablet. The latter view is supported by the
fact that the lines separating the columns on the obverse of the tablet are
continued on the reverse. Also, if the tablet is to list all ratios leading
to angles between $45$ and $30$ degrees, then it should contain an extra line
since the sixteenth ratio 16/9 yields the triplet (175, 288, 337), with an
angle slightly greater than $31^\circ$. All of this suggests that the size of
the tablet may be related to the purpose behind it, and possibly to the
method employed in its construction.  In the remainder of this section, we
will closely examine the three main previously proposed procedures for the
selection of $p$ and $q$, weigh the pros and cons of each procedure and,
based on the conclusion reached, select the most likely procedure to account
for the size of the tablet in a way which is consistent with the words and
numbers inscribed on it.

\paragraph{Procedure 1.} In this procedure, first introduced by Price, the
regular numbers $p$ and $q$ are chosen so that $1 < q < 60$ and $1 < p/q <
\alpha_0$ [\citealt{Price}]. Since 54 is the largest regular integer less
than 60 and since 128 is the largest regular integer less than 54 times
$\alpha_0$, the conditions on $p$ and $q$ can be restated as
\begin{equation}
1 < p \leq 128 \quad \text{and} \quad 1 < q \leq 54. \label{Condpq}
\end{equation}
This leads to a set of 38 admissible ratios (see Column I of Table
\ref{Cont}), the first 15 of which, when written in descending order, match
exactly with those needed to produce the Plimpton tablet. From this set, only
the fourth ratio $125/54$ cannot be written as $p\rcp q$, where $p$ and $q$
are in the standard table of reciprocals. The simplicity of the procedure
plus the fact that it succeeds in producing the correct ratios make it doubly
appealing, though not without some shortcomings. One such shortcoming is that
for $p$ and $q$ satisfying (\ref{Condpq}), we get a total of 234 distinct
ratios and there is no telling of how the 38 admissible ratios were selected
and sorted unless all 234 ratios were written as sexagesimal
numbers.\footnote{Since 54 is the largest regular numbers less than 60, taking
  $q=60$ generates only one extra ratio, namely 1/30, but this does not
  change the set of admissible ratios.} More importantly, it is not explained
why $q$ has to be less than 60, especially since, according to Price, entries
in the first column of the preserved tablet are calculated by squaring the
result of the division of $p^2+q^2$ by $2pq$. In other words, if both $\rcp
p$ and $\rcp q$ are needed in the calculation of the table, then there is no
good reason for $p$ and $q$ to have different upper bounds. This problem
persists even if the $r$-method is used, since for every $r=p \rcp q$ we must
also compute $\rcp r = q \rcp p$.

\paragraph{Procedure 2.} In this procedure the ratios are chosen so that both
$p$ and $q$ are regular integers less than or equal to 125. This leads to 47
admissible ratios, of which the line corresponding to the twelfth ratio
$125/64=1;57\ 11\ 15$ is missing from the tablet, see Table \ref{Extra}. The
argument proposed by Joyce that this ratio may have been dismissed because it
leads to large values of $w$ and $d$ does not carry much weight since the
ensuing triplet (11529, 16000, 19721) is comparable to the triplet (12709,
13500, 18541) produced by the fourth ratio $125/54$ [\citealt{Joyce}]. This,
however, does not completely rule out the possibility that the ratio was
unintentionally overlooked. But even if that is the case, we still have to
address the unanswered question: Why was 125 chosen as the upper bound for
$p$ and $q$? A possible answer to this crucial question is offered by the
following quote from Neugebauer: `\dots The only apparent exception is
$p=2;05$ but this number is again well known as the canonical example for the
computation of reciprocals beyond the standard table.' [\citealt[p.
39]{Neugebauer2}] As in Procedure 1, there is no mentioning of how the 47
admissible ratios were sorted out from a total of 303 possible ratios.

\paragraph{Procedure 3.} This procedure was first proposed by Bruins and
later adopted and improved by Robson [\citealt{Bruins1}; \citealt{Robson2}].
According to Robson, the ratios were chosen so that neither $r = p \rcp q$
nor $\rcp r = q\rcp p$ has more than four sexagesimal places, with the total
number of places in the pair not exceeding seven. The number of ratios
satisfying these conditions is 18, of which only 15 found their way to the
Plimpton tablet. The decimal values of the three omitted ratios are 288/125,
135/64 and 125/64, and the lines corresponding to them are shown in Table
\ref{ExtRob} (the letter `a' is appended to the line number $n$ to indicate
that the line at hand should be inserted between line $n$ and the line that
follows). 
\begin{table}[ht]
\centering
\setlength{\extrarowheight}{-2pt}
\addtolength{\tabcolsep}{5pt}
\begin{tabular}{lrlrrr}
\hline \\[-16pt] 
       $r$&     $l$&               $d^2/l^2$&      $w$&      $d$& $n$ \\[3pt]
\hline \\[-16pt] 
2;18 14 24& 20 00 00&    1;52 27 06 59 24 09& 18 41 59& 27 22 49& 4a \\
2;06 33 45&  4 48 00& 1;40 06 47 17 32 36 15&  3 55 29&  6 12 01& 8a \\   
1;57 11 15&  4 26 40& 1;31 09 09 25 42 02 15&  3 12 09&  5 28 41& 11a \\
\hline
\end{tabular}
\caption{The three lines corresponding to the ratios 288/125, 135/64 and
  125/64.}  
\label{ExtRob}  
\end{table}
The main problem with this procedure is that it is not easily justifiable why
the maximum number of sexagesimal places in the pair $r$ and $\rcp r$ has to
be seven. We believe that this is a consequence of rather than the criteria
for choosing $r$ and $\rcp r$.  Our view is supported by the fact that apart
from the fourth ratio 125/54, every other ratio can be expressed as $p\rcp
q$, where $p$ and $q$ belong to the standard table of reciprocals. But even
if we accept that the pair $r$ and $\rcp r$ should not have more than seven
sexagesimal digits, there are still some issues with this procedure that must
be addressed.

First, Robson assumes that because they yield what she calls \emph{nice}
(small) Pythagorean triplets, the first and fifteenth ratios were preselected
by the scribe as upper and lower bounds for the remaining ratios. But if that
is the case, then certainly the ratio $5/3 = 1;40$, whose reciprocal is 0;36,
makes a better lower bound (see line 18 of Table~\ref{Cont}). This is more so
since it generates the \emph{nicer} triplet (15, 8, 17). Second, she talks
about the ancient Babylonians being totally oblivious to the notion of a
complete table, which is hard to believe in light of the fact that every line
corresponding to a standard ratio lying between the largest (top) and
smallest (bottom) ratios is included in the tablet.  Third, she excludes line
4a on the ground that its short side and diagonal are \emph{half a place} too
long, meaning that they contain tens in the leftmost sexagesimal place. But
if the length of the sides is important then so should be the size of the
ratio $d^2/l^2$, where a quick glance at Column I of the tablet shows that
the value of $d^2/l^2$ in line~10 is two full places longer than its value in
line 4a.  In addition, she argues that since the long side is two sexagesimal
places in every line of the table, line 11a was left out because its longer
side 4 26 40 has three sexagesimal places.\footnote{By allowing both $r$ and
  $\rcp r$ to have up to four sexagesimal digits, Robson obtains three more
  pairs which are then dismissed using similar reasoning.} While this is true
if we exclude the terminating zero in the longer side, the scribe must be
fully aware of this zero, especially since the longer (uninscribed) side in
lines 2, 5 and 15 does not end with zero. Finally, she dismisses line 8a
because its generating ratio is not easily derivable from the standard table
of reciprocals using attested OB techniques. The same argument is given as
another reason for dismissing line 4a. We think, however, that the first two
procedures provide a better explanation of why the two lines are missing from
the tablet.

When all three procedures are considered, Procedure 1 is more likely to have
been used by the author of the Plimpton tablet because it produces the
complete set of ratios based on a minimal number of tenable assumptions. As
for the other two procedures, the weakness of the second stems from its
inability to account for the missing ratio 125/64; while the third is plagued
by the many contentious rules proposed by Robson, who has otherwise done an
excellent job in explaining the Plimpton tablet and in putting it in its
proper context.

\section{A New Method for Reconstructing the Table}
In this section, we will show how the generating ratios $p/q$ can be obtained
by a new selection procedure, which is not only easy to implement but is also
consistent with extant OB mathematics. In addition, we will demonstrate how
upper bounds for $p$ and $q$ like those given by Price may be reached
depending solely on the approximation of root two being used.

\begin{figure}[ht]
\centering
\psset{xunit=0.25}
\psset{yunit=2.0}
\begin{pspicture}(0,-0.25)(39,3)
\psaxes[Dx=5, dx=5, Dy=0.5, axesstyle=frame,tickstyle=bottom](0,0)(39,2.75)
\listplot[plotstyle=dots]{
1 2.4
2 2.37
3 2.344
4 2.315
5 2.25
6 2.222
7 2.16
8 2.133
9 2.083
10 2.025
11 2.0
12 1.92
13 1.875
14 1.852
15 1.8
}
\listplot[plotstyle=dots]{
1 1.983
2 1.949
3 1.919
4 1.886
5 1.815
6 1.785
7 1.72
8 1.693
9 1.643
10 1.586
11 1.562
12 1.489
13 1.45
14 1.43
15 1.387
}
\listplot[plotstyle=dots]{
1 0.9917
2 0.9742
3 0.9585
4 0.9414
5 0.9028
6 0.8861
7 0.8485
8 0.8323
9 0.8017
10 0.7656
11 0.75
12 0.6996
13 0.6708
14 0.6559
15 0.6222
}
\listplot[plotstyle=dots,dotsize=2pt]{
16 1.778
17 1.687
18 1.667
19 1.62
20 1.6
21 1.562
22 1.5
23 1.481
24 1.44
25 1.422
26 1.406
27 1.389
28 1.35
29 1.333
30 1.28
31 1.25
}
\listplot[plotstyle=dots,dotsize=2pt]{
16 1.369
17 1.3
18 1.284
19 1.251
20 1.238
21 1.213
22 1.174
23 1.163
24 1.139
25 1.129
26 1.121
27 1.112
28 1.093
29 1.085
30 1.062
31 1.051
}
\listplot[plotstyle=dots,dotsize=2pt]{
16 0.6076
17 0.5475
18 0.5333
19 0.5014
20 0.4875
21 0.4612
22 0.4167
23 0.4032
24 0.3728
25 0.3595
26 0.3476
27 0.3344
28 0.3046
29 0.2917
30 0.2494
31 0.225
}
\listplot[plotstyle=dots,dotsize=1pt]{
32 1.2
33 1.185
34 1.125
35 1.111
36 1.08
37 1.067
38 1.042
}
\listplot[plotstyle=dots,dotsize=1pt]{
32 1.034
33 1.029
34 1.014
35 1.011
36 1.006
37 1.004
38 1.002
}
\listplot[plotstyle=dots,dotsize=1pt]{
32 0.1833
33 0.1707
34 0.1181
35 0.1056
36 0.07704
37 0.06458
38 0.04083
}
\rput(5,2.5){$p/q$}
\rput(5,1.55){$d^2/l^2$}
\rput(5,0.7){$w/l$}
\end{pspicture}
\caption{The ratios $p/q$, $d^2/l^2$ and $w/l$ for the complete table.}
\label{Ratio} 
\end{figure}
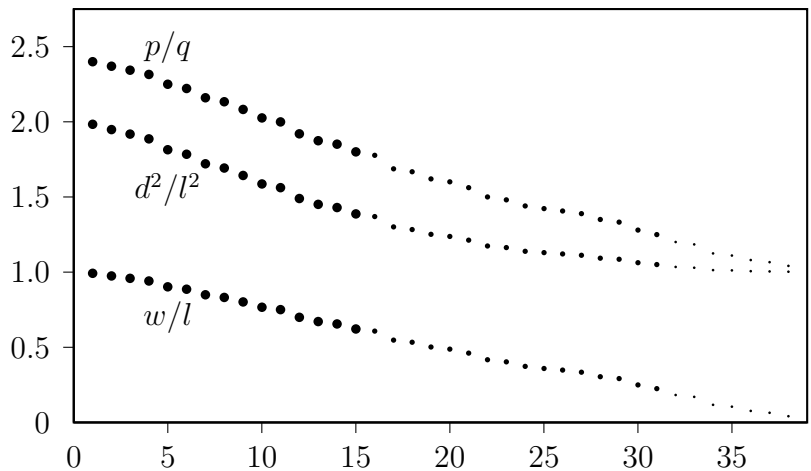

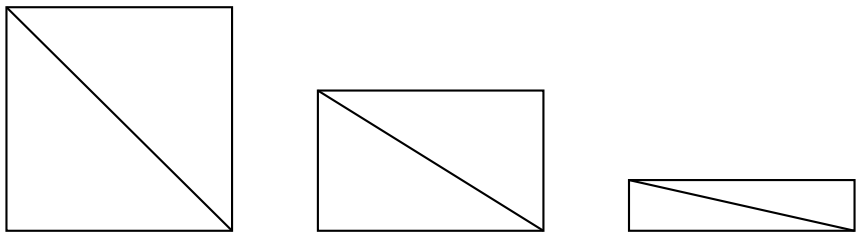
\begin{figure}[ht]
\centering
\psset{xunit=0.025}
\psset{yunit=0.025}
\begin{pspicture}(0,0)(120,120)
\pspolygon(0,0)(0,119)(120,119)(120,0)
\psline(0,119)(120,0)
\end{pspicture}
\hskip 1cm 
\psset{xunit=2.666}
\psset{yunit=2.666}
\begin{pspicture}(0,0)(45,45)
\pspolygon(0,0)(0,28)(45,28)(45,0)
\psline(0,28)(45,0)
\end{pspicture}
\hskip 1cm 
\psset{xunit=1.125}
\psset{yunit=1.125}
\begin{pspicture}(0,0)(40,40)
\pspolygon(0,0)(0,9)(40,9)(40,0)
\psline(0,9)(40,0)
\end{pspicture}
\caption{The rectangles corresponding to rows 1, 15 and 31 of the complete
  table.} \label{Full}
\end{figure}

We begin by noticing that the ratios $p/q$ and $d^2/l^2$ in columns one and
three of Table \ref{Cont} decrease as we move down the table. The same is
true for the ratio $w/l$, obtained by dividing an entry in Column IV by the
respective entry in Column II.  In Figure \ref{Ratio}, the three dotted
graphs (from top to bottom) represent $p/q$, $d^2/l^2$ and $w/l$: The large
dots correspond to lines $1$ to $15$, found on the obverse of Plimpton 322;
the midsize dots correspond to lines $16$ to $31$; and the small dots
correspond to lines $32$ to $38$.
The reason we grouped lines 32 to 38 together is that the ratio $l/w$ in
these lines is greater than 5, while for the first 31 lines we have
\begin{equation}
1 < l/w < 5. \label{Rev}
\end{equation}
More precisely, $l/w$ does not exceed 4 until line 30, and even for line 31
it is still less than half way between 4 and 5, 40/9 to be exact. So it seems
arguable that the table should end at line 31 since this leaves 16 entries
for the reverse of the tablet, as opposed to the 15 lines (plus headings)
inscribed on the obverse of the tablet. The argument is strengthened by the
fact that it is extremely difficult to find right triangles with $l/w > 5$ in
extant OB mathematics. On the other hand, there are many OB text problems
having right triangles (rectangles) whose lengths to widths are around 4 to 1
[\citealt{Robson1}]. In Figure~\ref{Full}, the rectangles corresponding to
lines 1, 15 and 31 are drawn. Observe that in addition to having their
lengths to widths satisfy~(\ref{Rev}), these triangles look like what a
teacher (whether ancient or modern) trying to illustrate Pythagoras theorem
would draw in front of a group of students.

\begin{figure}[ht]
\centering
\psset{xunit=0.01}
\psset{yunit=0.01}
\begin{pspicture}(0,0)(1200,100)
\pspolygon(0,0)(0,49)(1200,49)(1200,0)
\psline(0,49)(1200,0)
\end{pspicture}
\caption{The last (38-th) rectangle.} \label{Last}
\end{figure}
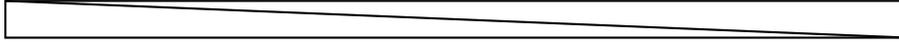

If an inequality similar to (\ref{Rev}) is to hold for all entries in the
tablet, then one would have to disagree with Price that the reverse of the
tablet should contain the $23$ uninscribed lines. This view is supported by
the fact that even if we assume that the reverse of the tablet does not
contain the two lines occupied by the heading on the obverse, a maximum of 17
or 18 lines could fit on the reverse.  In the case of 18 lines, we get a
total of 33 triangles (rectangles) with $1 < l/w < 6$. In contrast, the last
(38-th) rectangle, shown in Figure \ref{Last}, has a ratio $l/w$ greater than
$24$.  Therefore, if the Plimpton tablet is to be a list of practical
Pythagorean triplets then it should not contain such a triangle.
However if the tablet is thought of as an ancient piece of number theory then
there is no reason why the process should not be continued until all triplets
have been found. We believe that the size of the tablet and the nature of
Babylonian mathematics speak in favor the former point of view.  This does
not mean that the scribe was unaware of how the process can be carried out to
its fullest. As a rule of thumb, the scribe had to perform a balancing act
between how close is the first generating ratio to $\alpha_0$ (regular
numbers can get arbitrarily close to $\alpha_0$) and how large is the
corresponding Pythagorean triplet, keeping in mind that the ratio $l/w$ of
the resulting triangle should not be permitted to increase boundlessly.

For instance, suppose that the scribe considered all regular integers $p$ and
$q$ that are less than or equal $60^2$. Then the only ratio $p/q$ that falls
between the first ratio $12/5 = 2;24$ and $\alpha_0$ is $3125/1296 =
\text{2;24 40 33 20}$, which shows that $2;24$ is an exceptionally good
choice for the first ratio.\footnote{Even if we allow $p$ and $q$ to be as
  large as $60^3$, only one more ratio lying between 2;24 and $\alpha_0$ can
  be found, namely $19683/8192 = 2;24$ 09 45 21 05 37 30.} Consequently,
whether the scribe was searching for a regular number greater than 2;24 or
just considered regular numbers from a standard (or even extended) table of
reciprocals, the first ratio would be the same provided that the length of
the generated line should not exceed the width of the tablet. In particular,
taking $3125/1296$ as the first ratio, the first triplet becomes (8086009,
8100000, 11445241), with the ensuing line shown in Table \ref{Extra1}. It is
doubtful that the scribe had performed the tedious calculations necessary to
produce such a line; but even if he did, the sheer size of the numbers
involved gave him a compelling reason to reject it.
\begin{table}[ht]
\centering
\small
\begin{tabular}{lllllll}
\hline
 $p$&  $q$&     $l$&              $d^2/l^2$&     $w$&     $d$& $n$ \\
\hline
52 05& 21 36& 37 30 00 00& 1;59 47 34 27 27 58 38 07 21 36& 37 26 06 49& 52
 59 14 01&  1 \\ 
\hline
\end{tabular}
\caption{The line corresponding to $p\rcp q = \text{2;24 40 33 20}$.}
\label{Extra1}  
\end{table}
All of this forces us to be overly cautious as we try to rediscover how the
generating ratios may have been chosen by a ancient scribe whose mathematical
tools and interests are in many ways alien to ours. It should also be
emphasized that when we speak of `the scribe' we mean by that the group of
professionals that has developed these tools, possibly over a period longer
than the lifespan of any individual member of the group.

\begin{figure}[ht]
\centering
\psset{xunit=0.75}
\psset{yunit=45}
\begin{pspicture}(0,-0.01)(16,0.11)
\psaxes[Dx=1, dx=1, Dy=0.05, axesstyle=frame,tickstyle=top](0,0)(16,0.11)
\listplot[plotstyle=polygon,fillstyle=solid,fillcolor=lightgray]{
1 0
1 0.03424
2 0.03424
2 0
2 0.03036
3 0.03036
3 0
3 0.1038
4 0.1038
4 0
4 0.02981
5 0.02981
5 0
5 0.06521
6 0.06521
6 0
6 0.02727
7 0.02727
7 0
7 0.05004
8 0.05004
8 0
8 0.05655
9 0.05655
9 0
9 0.02362
10 0.02362
10 0
10 0.07308
11 0.07308
11 0
11 0.0394
12 0.0394
12 0
12 0.01978
13 0.01978
13 0
13 0.04308
14 0.04308
14 0
14 0.01794
15 0.01794
15 0
}
\end{pspicture}
\caption{A chart of the difference between consecutive values of $d^2/l^2$.}
\label{Diff} 
\end{figure}
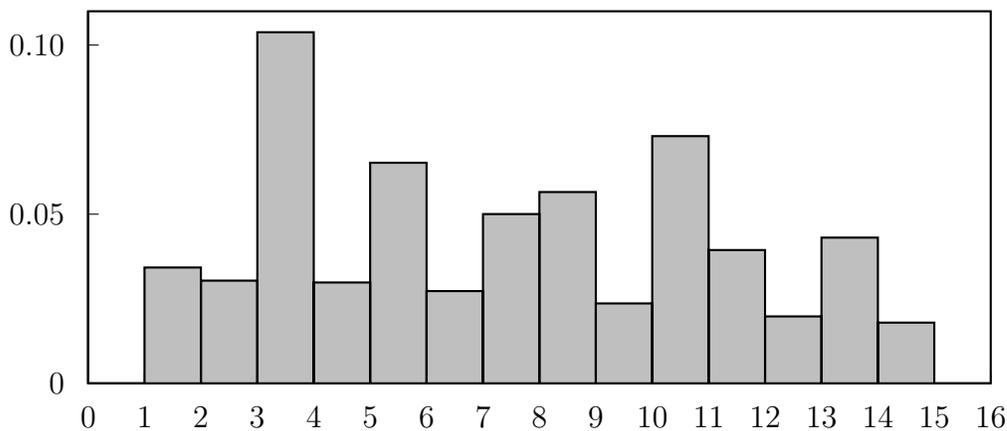

\paragraph{Procedure 4.} In this procedure, the generating ratios are those
of the form $p\rcp q$, where both $p$ and $q$ are taken from the standard
table of reciprocals (regular numbers up to 81).  This leads to 40 admissible
ratios, of which 81/64, 75/64 and 81/80 have $q > 60$. The three lines
produced by these ratios are listed in Table \ref{ExtReg}. As for the
remaining 37 lines, they are exactly those in Table~\ref{Cont}, with line 4
deleted.
\begin{table}[ht]
\centering
\setlength{\extrarowheight}{-2pt}
\addtolength{\tabcolsep}{5pt}
\begin{tabular}{lllrll}
\hline \\[-16pt] 
       $r$&     $l$&                 $d^2/l^2$&     $w$&      $d$& $n$ \\[3pt]
\hline \\[-16pt] 
1;15 56 15& 2 52 48& 1;03 23 29 29 33 54 01 40&   41 05&  2 57 37& 30 \\
1;10 18 45& 2 40 00&    1;01 31 19 18 53 26 15&   25 29&  2 42 01& 34 \\   
   1;00 45& 3 36 00& 1;00 00 33 20 04 37 46 40&    2 41&  3 36 01& 40 \\
\hline
\end{tabular}
\caption{The three extra lines corresponding to the ratios 81/64, 75/64 and
  81/80.}  
\label{ExtReg}  
\end{table}
Assuming that line 4 was not added to the tablet by mistake, it is plausible
that it was inserted between line 3 and line 5 to reduce the conspicuously
large difference between the third and fifth ratio. The motivation behind this
is that among the 40 admissible ratios, the largest difference between
successive ratios is 0;05~37~30 and the largest difference between successive
values of $d^2/l^2$ is 0;06 13 39 35 33 45, both of which occurring between
the third and fourth lines.  In Figure~\ref{Diff}, we draw the difference
between successive values of $d^2/l^2$ for the first fifteen ratios, where it
is obvious that the greatest difference is the one between line~3 and line~4.
In fact, if we restrict ourselves to the first 15 ratios, then 0;05~37~30 is
the only difference in $r$ exceeding $\rcp{12} = 0;05$, while
0;06~13~39~35~33~45 is the only difference in $d^2/l^2$ exceeding $\rcp{10} =
0;06$. Given the importance of the numbers 10 and 12 in the sexagesimal
number system, these facts can hardly be ignored.\footnote{The only other
  difference in $r$ exceeding $\rcp{12}$ is 0;05 25, taking place between the
  fifteenth and sixteenth ratios; while no other difference in $d^2/l^2$
  exceeds $\rcp{10}$.}  Moreover, the first 15 ratios are exactly those
leading to $1 > w/d > 0;30$, that is they cover every rectangle whose width
is greater than half its diagonal.\footnote{If lines 16 to 31 were inscribed
  on the reverse of the tablet, then they would be exactly those lines
  satisfying the inequality $0;30 > w/d > 0;12$.}  Since it is highly likely
that the tablet is a copy of an older original, one can see how the original
might have been produced by precisely these ratios, but line 4 was added to
the tablet either inadvertently by an inexperienced scribe or overtime by
someone who has noticed the unusually large gap between the third and fourth
ratios.

We still have to answer how the admissible ratios were sorted out. We have
seen that the standard table of reciprocals is comprised of the 30 regular
numbers less than or equal to 81. It follows that there is a total of 900
regular numbers of the form $r=p\rcp q$, where $p$ and $q$ are taken from the
standard table of reciprocals. Of these 900 ratios, 237 are distinct, and so
it is not an easy task to order them.\footnote{If we allow $p$ and $q$ to be
  1, then we get a total of 961 fractions, of which 257 are distinct.} The
simplest way would be to write down the sexagesimal representation of each
fraction and then sort them in ascending or descending order. But the number
of possible ratios can be greatly reduced if we are only interested in ratios
between some given bounds. For example, there are only 49 values of $r$
satisfying the condition $1 < r < 3$. An upper bound of $3$ is reasonable
because $3$ is the smallest integer greater than $\alpha_0$. In this case,
finding the admissible ratios amounts to choosing $p$ between $q$ and $3q$, a
condition that can be easily checked.

We believe that the above procedure provides a simple and direct way of
producing and sorting the generating ratios. Moreover, the sorting part of
the procedure can be modified so that it can be applied to other procedures.
In particular, Procedure 1 could become much more plausible provided that we
can satisfactorily explain why 60 (apart from being the base number) was
taken as an upper bound for $q$. Indeed, it can be shown that upper bounds
for $p$ and $q$ similar to those in (\ref{Condpq}) can be reached in a number
of ways, depending only on the approximation of $\sqrt2$ used. To begin with,
suppose that the scribe used 1;30 as an estimate of $\sqrt2$. Then his
approximation for $\alpha_0 = 1 + \sqrt2$ would be $2;30$. But in the
Babylonian number system $2;30$ is written as $2 \ 30$, which in addition to
our 5/2 could be read as 150 or even 150/60. The last form points us in the
direction of how the scribe may have obtained $2\ 30$ and $1 \ 00$ as upper
bounds for $p$ and $q$, leading to the same 38 ratios obtained using
(\ref{Condpq}).  Now the advantage of 2~30 as an upper bound is that it is
the smallest regular integer $p$ such that $p$ times $\rcp{60}$ is greater
than $\alpha_0$. It follows that for $q \leq 60$, no new admissible ratios
will be produced for values of $p$ larger than 2 30. Moreover, to determine
whether the ratio $p \rcp q$ is less than $2;30$, all that the scribe had to
do is check if $2p < 5q$. Alternatively, had the scribe used $1;25$ as an
approximation of $\sqrt2$, his estimate of $\alpha_0$ would be $2;25 =
29/12$. Again, he could think of this as $2 \ 25$ divided by 1 00, which in
turn could be taken as upper bounds for $p$ and $q$. As in the previous case,
the same admissible ratios are generated, but checking for admissibility in
this case is not as simple. This may not be a disadvantage since other tests
for admissibility, such as (\ref{Cond}), may have been used.

In addition to the two approximations of $\sqrt2$ used above, the scribe may
have used a third approximation, albeit indirectly. Earlier in this section,
we have seen that even for $p$ and $q$ as large as $60^2$, the only regular
number $p\rcp q$ lying between 2;24 and $\alpha_0$ is 2;24 40 33 20. Since
every ratio used in the tablet has a maximum of four sexagesimal digits and
since no such ratio exists between 2;24 and $\alpha_0$, the scribe may have
been prompted from the outset to consider 2;24 as the largest admissible
ratio. In fact, even if the scribe was oblivious to all of this, the same
conclusion could be reached by taking the first two sexagesimal digits of
$\alpha_0$, provided that $\sqrt2$ is approximated by $1;24 \ 51 \ 10$. But
if we only consider ratios $p \rcp q$ less than or equal to 2;24 and insist
that $q \leq 60$, then $p$ should not exceed $144$.  As in the previous two
cases, we get the 38 ratios listed in Table~\ref{Cont}.

\section{Explaining the Errors} \label{ErrorSec} The Plimpton tablet contains
a number of errors that when carefully analyzed may give us a clearer
understanding of how the numbers on the tablet were generated. The apparent
errors can be divided into two categories: Typographical errors and
computational errors.\footnote{Neugebauer considered only two errors; Friberg
  added two more; and Robson added another one.} The typographical errors,
shown in Table~\ref{ErrorS}, can be easily explained.
\begin{table}[ht]
\centering
\setlength{\extrarowheight}{-2pt}
\addtolength{\tabcolsep}{5pt}
\begin{tabular}{rrrrr}
\hline
Error& Line& Column&  Inscribed number& Correct number \\
\hline
    1&    2&      I&    58 14 {\bf 56} 15&  58 14 {\bf 50 06} 15   \\
    2&    8&      I& 41 33 {\bf 59} 03 45& 41 33 {\bf 45 14} 03 45 \\
    3&    9&     II&           {\bf 9} 01&              {\bf 8} 01 \\
    4&   13&      I&             27 \ 03 45&       27 {\bf 00} 03 45 \\
\hline
\end{tabular}
\caption{The typographical errors in Plimpton 322.} \label{ErrorS} 
\end{table}
Looking at the first error, it is obvious that the scribe carelessly wrote
the symbol for 6 a bit too close to that of 50, and thus the correct number
can be obtained by simply inserting a little space between the two symbols.
To undo the second error, the sexagesimal digit 59 should be written as 45
followed by 14, and not as the sum of the two digits.\footnote{Friberg offers
  a reasonable explanation of how the error may have occurred if the
  $r$-method was used [\citealt{Friberg1}].} The third error amounts to 8
being miscopied as 9, where it is quite easy to make such a mistake due to
the similarity between the two symbols (9 has just one more wedge than 8).
As for the fourth error, it can be easily discarded if the space between 27
and 3 is transliterated as a zero. At any rate, this is hardly an error since
we know that OB scribes did not consistently use a blank space to represent
zero.

Having dealt with the typographical errors, the computational errors, listed in
Table \ref{ErrorC}, require a deeper understanding of OB mathematics. 
To begin with, we have seen that there is a definite error in the last line
of the tablet, but scholars are divided on whether the erroneous number
should be in the second or third column. More precisely, either the entry in
Column~III should be 1~46 rather than 53, or the entry in Column II should be
28 instead of 56. This suggests that at some point in the calculation, a
multiplication by 2 or $\rcp2$ should have been applied to both entries, but
the operation was only performed on one. But if the $pq$-method is employed,
then no doubling or halving is needed in calculating $p^2 - q^2$ (Column II)
or $p^2 + q^2$ (Column~III). This led proponents of the method to propose
that in order to obtain a primitive pair, the scribe intended to halve both
numbers but forgot to do it for $p^2-q^2$. On the other hand, if the
$r$-method was used, then the process of eliminating the common factors of
$x=0$;37~20 and $y=1$;10 40 should look something like this:
\begin{table}[t]
\centering
\setlength{\extrarowheight}{-2pt}
\addtolength{\tabcolsep}{5pt}
\begin{tabular}{rrrrr}
\hline
Error& Line& Column&  Inscribed number& Correct number \\
\hline
    1&    2&    III&           3 12 01&        1 20 25 \\
    2&   13&     II &          7 12 01&           2 41 \\
    3&   15&     II&                56&             28 \\
\hline
\end{tabular}
\caption{The computational errors in Plimpton 322.} \label{ErrorC} 
\end{table}

\begin{table*}[ht]
\centering
\setlength{\extrarowheight}{-2pt}
\begin{tabular}{rrlllll}
$3$&$\times$& 0;37 20 & \qquad 1;10 40& $\times$& $3$ \\
$30$&$\times$& 1;52 & \qquad 3;32& $\times$& $30$ \\
$30$&$\times$& 56& \qquad 1;46& $\times$& $30$ \\
&&28 00& \qquad 53 00&& 
\end{tabular}
\end{table*}
\vskip -0.5cm
\noindent
Since the elimination process is usually carried out on a supplementary
tablet, it is easily seen how 56 and 28 could be confused for each other as
the bottom numbers are transcribed from the supplementary tablet to the
Plimpton tablet. 

The above $r$-method explanation of the error is more consistent than the one
given by Robson, where the common factors of $x$ and $y$ are eliminated using
the multipliers 3, 5 and 3. The problem with Robson's explanation is that it
is not clear why the pair (1;52,~3;32) was multiplied by 5 rather than 30,
especially since in explaining the first computational error she multiplied
the pair (1;56 54 35, 2;47 32 05) by 12 and not by two.\footnote{Robson's
  justification is that the numbers in the first pair terminate in 2, while
  those in the second pair terminate in~5 [\citealt[pp.~192-193]{Robson2}].}
To be consistent, she should either multiply the latter pair by 2, or
preferably multiply the former pair by 15 so that, at each step, the product
of the greatest common factor and the multiplier is sixty. But if the second
multiplier is taken to be 15, then the process terminates after only two
steps, putting our choice of multipliers in question. This, however, can be
answered in one of two ways: Either the scribe failed to notice that 4 is the
greatest common divisor of 52 and 32, or he was aware of this but, since both
numbers end in 2, he went for the simpler operation of halving, which is
equivalent to multiplying by 30.  In the former case, the error probably
occurred at the second step when, rather than the single correct multiplier
of 15, the two multipliers 30 and 15 were respectively applied to the current
numbers 1;52 and 3;32, perhaps because noticing that 52 is a multiple of 4 is
not as obvious as noticing that 32 is a multiple of 4 [\citealt{Friberg1}].
But in either case, the ensuing explanation is more concise than the one
advanced by Robson.

Turning to the second computational error, we find that the inscribed number
7~12~01 is simply the square of the correct value 2 41. This, coupled with
the assumption that the entry in Column~III of line 15 should be twice the
inscribed number, led Gillings to conclude: `Thus to calculate the numbers of
the tablet, the doubling of numbers, and the recording of squares, from their
abundant table texts, must have been part of the procedure.'
[\citealt{Gillings2}] The doubling part of Gillings' conclusion is disputed
by the fact that our above explanation of the error in line 15 did not
require the doubling of numbers to be a necessary part of the procedure. In
addition, Gillings does not specify the stage of the procedure at which the
doubling and squaring should take place. If anything, these facts should
encourage us to search for new ways to explain how the square of the supposed
number made it to the tablet.  Although the exact way may never be known, we
can still make an educated guess. The first thing that comes to mind is that
the error may be due to a routine check that the scribe performed to make
sure that the numbers in Columns II and III along with the uninscribed side
form a Pythagorean triplet. To do so, the scribe first computes $w^2$ and
$d^2$ (after $w$ and $d$ were found by eliminating the common factors of $x$
and $y$) either directly or by consulting a table of squares. Then he
computes the difference between the two squares, and takes the square root of
the answer to get $l$.  At this stage, he is ready to transfer the results
from the rough to the clean tablet, but in the process of doing so he copied
$w^2$ instead of $w$.  This does not undermine our explanation of the error
in line 15 since transcribing the values of $w$ and $d$ on the clean tablet
should not be done until after the prescribed check has been performed. An
objection to this argument would be that when the $r$-method is used, then at
exactly the same step we reach $w$ and $d$, we also get $l$, which is equal
to the product of the multipliers used. A possible answer to this is that the
scribe cleared the common factors of $y$ only, and then found $w$ as the
square root of the difference between $d^2$ and $l^2$. But even if the scribe
knew how to find $l$ at the same time he found $w$ and $d$, he may still
choose to find it using the rule of right triangle, as we shall see in the
next section.

We are left with the first computational error, which is the most difficult
one to explain since there is no obvious relation between the inscribed
number 3 12 01 and the correct number 1~20~25. According to Neugebauer:
\begin{quote}
It seems to me that this error should be explicable as a direct consequence
of the formation of the numbers of the text. This should be the final test
for any hypothesis advanced to explain the underlying theory.   
\end{quote}
The above quote is taking from a note on page 50 of the 1951 edition of
Neugebauer's now classic book \emph{The Exact Sciences in Antiquity}. In
later editions of the same book, Neugebauer changed the note so that it
reflects R.  J.  Gillings attempt to resolve the error
[\citealt{Neugebauer2}].  According to Gillings, the error is due to the
accumulation of two mistakes made by the scribe [\citealt{Gillings1};
\citealt{Gillings4}]. First, in computing
\begin{equation}
d = p^2+q^2 = (p+q)^2 - 2pq, \label{MultSqre}
\end{equation}
the scribe accidently calculated $(p+q)^2+2pq$. For $p=1 \ 04$ and $q = 27$,
the calculated value would be
\[
d = 2 \ 18 \ 01 + 57 \ 36 = 3 \ 15 \ 37,
\]
but the scribe made the second mistake where he took $p = 1 \ 00$ instead of
1 04, obtaining $2pq = 54 \ 00$. Now adding 2 18 01 to 54 00, the inscribed
value is reached. A modified version of Gillings' conjecture was proposed by
Price, but Price considered the problem as unresolved. Later Gillings refuted
the equivalence between his method and that of Price [\citealt{Gillings4}].

The problem with Gillings' conjecture is twofold. First, the idea of
computing $p^2 + q^2$ as $(p+q)^2 - 2pq$ is doubtful since in this case $p=1
\ 04$ and $q=27$ are powers of $2$ and $3$ respectively, and so the scribe
can find $d$ by simply finding the sum of the easily computable squares of
$p$ and $q$. But even for arbitrary $p$ and $q$, the argument that the right
hand side of (\ref{MultSqre}) should be used to calculate $d$ begs the
question since it runs from finding the squares of $p$ and $q$ to finding the
square of the larger number $p+q$.  Moreover, for the Babylonians, squaring
was one of the basic mathematical operation, which is reflected not only in
the number but also in the scope of the tables of squares they left behind
[\citealt[45--52]{Friberg2}]. In fact, it is known that the Babylonians used
squares to find the product of two numbers by applying the formula
\begin{equation*}
pq = \frac14 \left[ (p + q)^2 - (p - q)^2 \right] \qquad \text{or} \qquad pq =
\frac{1}{2} \left[(p + q)^2 - p^2 - q^2\right], \label{SqreMutl}  
\end{equation*}
see [\citealt{Connor}]. This clearly supports our argument since in both
formulas the squares of $p$ and $q$ are used to calculate the product $pq$
and not the reverse, as proposed by Gillings. Second, the $pq$-method, and
hence Gillings' explanation of the error, involves the concept of
\emph{relatively prime} integers, which is not only unattested from the
historical and archaeological record but also runs contrary to the practical
nature of OB mathematics. So, in order to accept Gillings' conjecture, we
have to accept that a questionable method was used to generate the tablet;
that an unlikely rule was used to compute the sum of two squares; and that at
some point in the calculation subtraction was substituted for addition, and
then 1 04 was inexplicably taken as 1 00. In light of these facts, we think
that Gillings' explanation of the first computational error is improbable at
best.

The other often mentioned method for explaining the error in the second line
of the Plimpton tablet is based on the $r$-method introduced by E. M. Bruins
in 1949. The advantage of this method is that it employs attested OB
techniques, which made it the preferred method for authors like Friberg,
Schmidt, Robson and others.  Starting with the initial pair $x=0$;58~27~17~30
and $y=1$;23 46 02 30, Robson applied a slightly different
version of the method to obtain the numbers inscribed in columns two and
three as follows:

\begin{table*}[ht]
\centering
\setlength{\extrarowheight}{-2pt}
\begin{tabular}{rrrlrllll}
 2& $\times$& 0&\hskip-0.4cm;58 27 17 30& \qquad 1&\hskip-0.4cm;23 46 02 30&
 $\times$& 2  \\ 
12& $\times$&    1&\hskip-0.4cm;56 54 35& \qquad    2&\hskip-0.4cm;47 32 05&
$\times$& 12 \\ 
12& $\times$&      23&\hskip-0.4cm;22 55& \qquad      33&\hskip-0.4cm;30 25&
$\times$& 12 \\ 
12& $\times$&       4 40&\hskip-0.4cm;35& \qquad       6 42&\hskip-0.4cm;05&
$\times$& 12 \\ 
12& $\times$&         56 07&& \qquad       1 20 25&& $\times$& 12 \\
\hline
12& $\times$&    11 13 24&& \qquad        16 05 00&& $\times$& 12 \\
  &         &  2 14 40 48&& \qquad      3 13 00 00&&         & 
\end{tabular}
\end{table*}
\vskip -0.5cm
\noindent
Instead of stopping when the numbers 57 07 and 1 20 25 were reached, the
scribe, unaware that all common factors have been cleared out, carried out
the process two extra steps obtaining $w=2$ 14 40 48 and $d=3$ 13 00 00.
Realizing that he went too far, he looked back for the correct pair, but in
doing so he made two new mistakes: First, he took the value of $d$ from the
last step in place of the correct value obtained two steps earlier, and then
he sloppily wrote 3~12~01 instead of 3 13 [\citealt{Robson2}].  There is no
quarrel in accepting the assumption that 3 12 01 was written for 3 13 since
one will probably find similar mistakes somewhere in this
paper.\footnote{Despite the availability of modern computers, both
  typographical and computational errors are still being made by prominent
  authors in peer reviewed journals. We have seen that in his version of the
  complete table Price has made many calculation errors, some of which are
  similar to those in Plimpton~322 (in Column IV of line 15, 8 should be 9),
  while others are even more difficult to explain. A more serious mistake was
  committed by Friberg as he tried to explain the appearance of the number
  $M=2$~02~02~02~05~05~04 in the Sippar text Ist.S 428. Friberg erroneously
  observed, that the square root of $M$ is the same as the integral part of
  the square root of 2~02~02~02~02~02~02, while the two supposedly equal
  numbers are 1 25 34 {\bf08} and 1 25 34 {\bf07}
  [\citealt[290]{Friberg1}]. Even Robson has made a number of errors in her
  major work on the Plimpton tablet [\citealt{Robson2}]. Of these errors, two
  are extremely pertinent to us. The first error occurs in the bottom row of
  Table~6, where the numbers 1;48, 0;33 20 and 0;{\bf33} 20 are written under
  the headings $x$, $1/x$ and $(x-1/x)/2$. The first two entries are correct,
  but the third entry should have been 0;{\bf37} 20. The similarity between
  this and some of the errors on the Plimpton tablet is striking. But most
  ironic is the error made by Robson while trying to explain the first
  computational error. When multiplying 11 13 24 by 12, she wrote 2 14
  {\bf36} 48 instead of 2 14 {\bf40} 48. It seems that humans are still prone
  to the same mistakes they were prone to 4000 years ago.}  Even Gillings,
who staunchly opposed Bruins' explanation of the error, admitted that the
number 3 12 01 may well be 3 13 [\citealt{Gillings3}].  But the problem with
Bruins' method is that it is hard to see why the scribe would copy the
correct value for $w$ as opposed to the (wrong) value for $d$ that can only
be reached if two unnecessary steps have been performed. One explanation
would be if the multiplications were not carried out simultaneously. In such
case the error must have occurred as the scribe copied the numbers obtained
at the end of the simplification process, not taking into account that the
simplification of $y$ required six steps while that of $x$ required only four
steps.

The above procedure for explaining the first computational error can be
slightly modified so that it becomes much more plausible. Assuming that the
number 3 12 1 is a miscopy of 3 13, the values of $w$ and $d$ given by the
scribe can be reached in the following way:

\begin{table*}[ht]
\centering
\setlength{\extrarowheight}{-2pt}
\begin{tabular}{rrrlllll}
$\rcp{30}$& $\times$& 58 27 17 30& \qquad 1 23 46 02 30& $\times$&
$\rcp{30}$ \\ 
$\rcp 5$& $\times$& 1 56 54 35& \qquad 2 47 32 05& $\times$& $\rcp 5$ \\
$\rcp 5$& $\times$& 23 22 55& \qquad 33 30 25& $\times$& $\rcp 5$ \\
$\mathbf 5$& $\times$& 4 40 35& \qquad 6 42 05& $\times$& $\rcp 5$ \\
$\rcp 5$& $\times$& 23 22 55& \qquad  1 20 25& $\times$& $\rcp 5$ \\
$\rcp 5$& $\times$& 4 40 35& \qquad  16 05& $\times$& $\rcp 5$ \\
&&56 07& \qquad 3 13&& 
\end{tabular}
\vskip -0.5cm
\end{table*}
\noindent
In the first step the scribe correctly multiplied both $x$ and $y$ by
$\rcp{30}$, but in one of the following three steps (say step four) he must
have multiplied the left number by $5$ and the right number by $\rcp5$,
instead of multiplying both by $\rcp5$. The simplicity of the error and the
fact that this is the only case where the scribe needs four separate
multiplications in order to get rid of all common factors between $x$ and $y$
make the modified procedure much more attractive.

Another way to account for the error would be if the wrong value of $y$ was
used at the beginning of the procedure. For example, if $y=1$ 23 46 02 30 is
replaced by 3 20 01 02 30, then the inscribed number 3 12 01 would appear
opposite to 56 07. On the other hand, if $y$ is replaced by 3 21 02 30, then
the process terminates with 56 07 and 3 13. In both cases, the initial number
used is quite similar to the correct value of $y$. Since $y = \rcp2(r+\rcp
r)$, the error would have probably occurred as $r/2$ was added to $\rcp
r/2$. In particular, for $r=2;22\ 13\ 20$ the calculation of $y$ should look
something like this:

\begin{table*}[ht]
\centering
\setlength{\extrarowheight}{-2pt}
\begin{tabular}{l}
1;11 {\bf 06} 40 \\
0;12 {\bf 39} 22 30 \\ 
\hline
1;23 {\bf 46} 02 30
\end{tabular}
\end{table*}
\noindent
Now observe that if the (sexagesimal) digit 46 of the correct $y$ is ignored,
then the resemblance between 1 23 02 30 and 3 20 01 02 30 becomes more
pronounced when the two numbers are read out loud. It is even more so between
1 23 02 30 and 3 21 02 30. Moreover, 3 21 02 30 is equal to 12/5 times the
correct $y$, where 12/5 is the value of $r$ in the previous line.\footnote{If
  the calculations are carried out on a separate tablet, as is usually the
  case, then it is not difficult to see how $y=\rcp2(r+\rcp r)$ could have
  been multiplied by $12/5$, the value of $r$ in the previous
  calculation~of~$y$.} Applying Robson's procedure with $y=3\ 21\ 02\ 30$
yields: 

\begin{table*}[ht]
\centering
\begin{tabular}{rrrlrllll}
$2$&$\times$& 0&\hskip-0.4cm;58 27 17 30& \qquad 3&\hskip-0.4cm;21 02 30&
$\times$& $2$ \\ 
$12$&$\times$& 1&\hskip-0.4cm;56 54 35& \qquad 6&\hskip-0.4cm;42 05&
$\times$& $12$ \\ 
$12$&$\times$& 23&\hskip-0.4cm;22 55& \qquad 1 20&\hskip-0.4cm;25& $\times$&
$12$ \\ 
$12$&$\times$& 4 40&\hskip-0.4cm;35& \qquad 16 05&\hskip-0.4cm& $\times$&
$12$ \\ 
&&56 07&\hskip-0.4cm& \qquad 3 13 00&& 
\end{tabular}
\end{table*}

From the above discussion we see that except when we take $y=3\ 20\ 01\ 02\
30$, the corresponding number on the tablet should be 3 13 (or 3 13 followed
by one or two zeros) rather than the apparent 3 12 01. As to how this might
have happened, we offer three different explanations. First, the space
between 12 and 01 should be ignored as a scribal error, and consequently 3
12~01 should be read as 3 13. Second, being aware of the elusive zero(s) at
the end of 3 13, the scribe inserted a space (the OB symbol for zero) in the
wrong place and wrote 3 12 01 or what might be 3 12 00 01. Third, as the
scribe was copying the results on the tablet, he noticed that $56\ 07$ (the
width) is larger than $3\ 13$ (the diagonal). Realizing that this cannot be
true, he glanced over his calculations and hastily misread 3 13 as 3 12 01.

Having looked at the computational errors one by one, we now look at the
tablet as a whole to see if there is a possible relation between the
different errors. In Table \ref{New}, we list the lines for which the
rightmost digit of at least one of the (correct) values of $w$ and $d$
(Columns II and III of the tablet) is divisible by a regular number. 
\begin{table}[ht]
\centering
\setlength{\extrarowheight}{-2pt}
\addtolength{\tabcolsep}{3pt}
\begin{tabular}{llllll}
\hline 
        $x$&           $y$&   $l$&   $w$&     $d$& $n$\\
\hline
58 27 17 30& 1 23 46 02 30& 57 36& 56 07& 1 20 25&   2\\
      54 10&       1 20 50&  1 12&  1 05&    1 37&   5\\
         45&          1 15&  1 00&    45&    1 15&  11\\
      37 20&       1 10 40&    45&    28&      53&  15\\
\hline
\end{tabular}
\caption{The four lines for which $w$ or $d$ is divisible by a regular
number.}   \label{New}  
\end{table}
Observe that except for line 5, one or both of $w$ and $d$ deviate from the
expected answer. For lines 2 and 15, the errors would be easily explained if
the calculation of $w$ from $x$ was not done in parallel with the calculation
of $d$ from $y$. But even if the calculations were done side by side, it is
still not difficult to see how the numbers inscribed on the clean tablet
could be out of step with each other, especially since for these lines the
number of steps needed to clear out the regular numbers on the $x$ side is
different from the number of steps needed on the $y$ side. As for line 11,
the numbers inscribed in Columns II and III are those of $x$ and $y$,
respectively 45 and 1 15. Since in this case the value of $l$ is the base
number 1~00, the scribe was apparently satisfied with the non-reduced triplet
(45,~1~00,~1~15), which Melville called the favorite version of the primitive
triplet (3, 4, 5) [\citealt{Melville1}]. Moreover, since we are fairly sure
that the scribe knew of the equivalence between the two triplets, it not
unreasonable to think that he kept the non-reduced triplet because, unlike
the primitive triplet, its entries are comparable in magnitude to those of
other triplets.

\section{Purpose of the Tablet}
The presence of modern glue on the edge of the Plimpton tablet does not
completely rule out the possibility that the tablet was damaged while it was
still in preparation. In such case, the idea that the reverse was meant for
the remaining entries becomes much more plausible. Alternatively, the tablet
may have been abandoned because someone had noticed the errors before the
reverse had been inscribed. This goes hand in hand with the view that the
tablet is just a school exercise about finding right triangles with integral
sides, which is normally giving by an accomplished scribe to a group of
students who aspire to join the respected profession of their teacher. In
fact, the types of computational errors committed added to the similarity
between two of the errors suggest in a way that the tablet may have been
written by someone who has not fully mastered the techniques involved.

If the tablet is taken to be an exercise from a scribal school, then it is
not that hard to see the purpose behind such an exercise.  First, the
generation of the tablet involves many of the mathematical techniques used by
OB scribes: Cut-and-paste geometry, division as multiplication by a
reciprocal, squaring and taking the square root of a number, and so on.
Second, it is easy for the presiding scribe to check the work of his students
by comparing their results with a master copy that has the correct answers.
This will be more so if the tablet contains a complete list of triplets,
meaning all triplets corresponding to standard ratios between the first and
last ratio of the list, as is the case with the Plimpton tablet.  Third, the
numbers on the tablet may be related to the solution of another ancient
problem about right triangles and upright walls. True, there is no direct
evidence that the Babylonians used the rule of right triangle in erecting
walls, but there are extant Babylonian (and even Egyptian) problem texts in
which the rule is used to measure the length of a cane leaning against a wall
[\citealt{Melville1}].

\begin{figure}[ht]
\centering
\psset{xunit=0.2cm}
\psset{yunit=0.2cm}
\begin{pspicture}(-2,-0.75)(10,15)
\pspolygon[fillstyle=solid,fillcolor=lightgray](-2,0)(0,0)(0,15)(-2,15)
\psline[linewidth=3pt](0.2,12)(9,0)
\psline(-2,0)(12,0)
\rput(4.5,-1.33){$b$}
\rput(4.5,8.5){$l$}
\rput(1,6){$h$}
\rput(0.75,14){$d$}
\end{pspicture}
\caption{The cane against the wall problem of BM 34 568.}
\label{Cane} 
\end{figure}
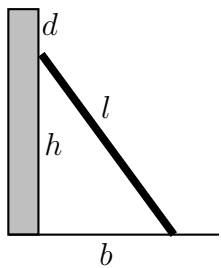

One such Babylonian cane-against-the-wall problem is found in the BM 34 568
tablet from the Seleucid period, roughly 300--100 BCE. The problem as stated
by Friberg reads like this:
\begin{quote}
  A cane is leaning against a wall. 3 cubits it has come down, 9 cubits it
  has gone out. How much is the cane, how much the wall? I do not know their
  numbers [\citealt{Friberg1}]. 
\end{quote} 
Following the statement of the problem, the length of the cane is then found
by calculating
\[
l = \frac{d^2+b^2}{2d}, 
\]
where for $d=3$ and $b=9$ we get $l=15$, see Figure \ref{Cane}.
Using modern notation, the solution is obtained by writing the equation
$h^2+b^2=l^2$ as
\[
(l-d)^2 + b^2 = l^2 \quad \text{or} \quad d^2 + b^2 = 2dl,
\]
and then solving for $l$. Finally, instead of calculating the height of the
wall using the formula $h=l-d = 15-3 = 12$, the scribe used the rule of right
triangle to first compute $h^2 = l^2-b^2 = 2 \ 24$, and then took the square
root of $2 \ 24$ to get $h= 12$. As mentioned in the previous section, this
gives more credibility to our explanation of the second computational error,
since the erroneous number inscribed on the tablet is simply the square of
the correct number.

Although the above cane-against-the-wall problem comes from period separated
by nearly fifteen centuries from the time of Plimpton 322, it is still
similar in style as well as content to tablets from the OB period. In fact, a
similar problem is found on the OB tablet BM~85196, where the triplet (18,
24, 30) provides the correct answer [\citealt{Melville1}].\footnote{Observe
  that in both problems, the resulting triplet is a multiple of the primitive
  triplet (3, 4, 5).} Moreover, Carlos Gon\c calves has recently shown that
the solution to the first problem of BM~34~568, which is about finding the
diagonal of a rectangle but is not solved using the rule of right triangle,
can be reduced to finding a pair of reciprocals either algebraically or
preferably using cut-and-paste geometry [\citealt{Goncalves}]. While the
usefulness of the cane-against-the-wall problem supports the view that the
Plimpton tablet can be thought of as a scribal school exercise with some
practical applications, the argument given by Gon\c calves clearly favors
Bruins' method of generating the numbers on the tablet.

\bibliographystyle{plainnat}

\end{document}